\newcommand{\here}{\noindent\hbox{$\blacktriangleright$}%
	    }
\newcommand{\there}{\noindent\hbox{$\blacktriangleleft$}%
	       }

\renewcommand{\here}{}
\renewcommand{\there}{}
\documentclass [12pt] {amsart}
\usepackage{latexsym}
\usepackage{amssymb}
\usepackage[latin1]{inputenc}
\usepackage{a4}
 \usepackage{amsthm}
\usepackage{mathrsfs}
 \newcommand{\cal}{\mathscr}
\newcommand{\halfopen}[2]{ [ #1 , % ]( 
                             #2 ) }

\newcommand{\anzahl}[2]{Z_{#1}({#2})}

\theoremstyle{remark}
\newtheorem{Def}{Definition}[section]

\newtheorem{Fact}[Def]{Fact}
\newtheorem{Rem}[Def]{Remark}

\theoremstyle{plain}
\newtheorem{Theo}[Def]{Theorem}
\newtheorem{Theorem}[Def]{Theorem}
\newtheorem{Prop}[Def]{Proposition}

\newtheorem{Lemma}[Def]{Lemma}

\newcommand{\N}{{\mathbb N}}
\newcommand{\Z}{{\mathbb Z}}

\newcommand{\R}{{\mathbb R}}

\newcommand{\xf}{{\bf x}}
\newcommand{\yf}{{\bf y}}

\newcommand{\nf}{{\bf n}}

\newcommand{\If}{{\bf I}}
\newcommand{\mf}{{\bf m}}
\newcommand{\qf}{{\bf q}}

\newcommand{\tm}{\subseteq}
\newcommand{\xn}{(x_n)_{n \in \N}}
\newcommand{\yn}{(y_n)_{n \in \N}}

\newcommand{\Ij}{(I_j)_{j \in \N}}

\newcommand{\nk}{(n_k)_{k \in \N}}

\newcommand{\supp}{{\mbox{supp}}}

\newcommand{\MM}{{\cal M}}
\newcommand{\CC}{{\cal C}}
\newcommand{\PP}{{\cal P}}

\newcommand{\limn}{\lim_{n \to \infty}}
\newcommand{\e}{\varepsilon}
\newcommand{\ve}{\varepsilon}
\newcommand{\boundary}{\partial}

\newcommand{\forces}{\Vdash}

\title [Baire results on the distribution of subsequences]
  {Further Baire results\\ on the distribution of subsequences}

\date{August 2007}

\author{Martin Goldstern}
\author{J{\"o}rg Schmeling}
\author{Reinhard Winkler}
\email{martin.goldstern@tuwien.ac.at}
\email{reinhard.winkler@tuwien.ac.at}
\email{joerg@maths.lth.se}

\address[Goldstern, Winkler]{\endgraf
Institut f{\"u}r Diskrete Mathematik und Geometrie\endgraf
Technische Universit\"{a}t Wien\endgraf
Wiedner Hauptstra{\ss}e 8-10\endgraf
1040 Wien\\
Austria}

\address[Schmeling]{\endgraf
Center for Mathematical Sciences\endgraf
LTH, P.O.-Box 118 \endgraf
SE-22100 Lund \\
Sweden}

 \makeatletter
 \def\@secnumfont{\bfseries}

\begin{document}

\begin{abstract}
This paper presents results about the distribution
of subsequences which are typical in the sense of Baire categories.

The first main part is concerned with sequences of the 
type~$x_k=n_k\alpha$, $n_1 < n_2 < n_3 < \cdots$, mod 1. Improving
a result of {\v S}al{\'a}t we show that, if the quotients
$q_k = n_{k+1}/n_k$ satisfy $q_k \ge 1+\e$, then the set of all
$\alpha$ such that~$(x_k)$ is uniformly distributed is of first
Baire category, i.e.\ for generic~$\alpha$ we do not have
uniform distribution. Under the stronger assumption 
$\lim_{k \to \infty} q_k = \infty$ one even has maldistribution
for generic~$\alpha$, the strongest possible contrast to
uniform distribution. Nevertheless, growth conditions on the $n_k$
alone do not suffice to explain various interesting phenomena.
In particular, for individual sequences the situation maybe
quite diverse: For $n_k=2^k$ there is a set $M$ such that for
generic $\alpha$ the set of all limit measures of $(x_k)$ is
 exactly~$M$, while for $n_k=2^k+1$ such an $M$ does not exist.

For the rest of the paper we consider appropriately
defined Baire spaces~$S$ of subsequences. For a fixed well distributed
sequence~$(x_n)$ we show that there is a set~$M$ of measures such that
for generic $(n_k) \in S$ the set of limit measures of 
the subsequence~$(x_{n_k})$ is exactly~$M$. 
\end{abstract}



\thanks{RW thanks the Austrian Science Foundation FWF
for its support through Projects no.\ S8312 and S9612-N13.}

\thanks{A preprint of this paper 
is available at {\tt http://arxiv.org/math.NT/0407295 }.}
\dedicatory{Dedicated to Professor Robert F.~Tichy on the occasion of his 50th birthday}
\maketitle

\section{Introduction}\label{Sintro}

\subsection{Motivation}\label{SSmotiv}

This paper is a continuation of topological investigations
contained in \cite{GSW}.
Let~$X$ be a compact metric space and let
$\xf = \xn$, $\yf = \yn$ etc.\ denote sequences on~$X$, and
\[
A(\xf):= \bigcap_{n_0 \ge 1} \overline{\{x_n:\ n \ge n_0\}}
\]
the set of accumulation points of the sequence~$\xf$.
We are interested in subsequences of~$\xf$, therefore we write
$\nf = \nk$ for sequences $0 < n_1 < n_2 < \cdots$ of positive
integers and $\xf \nf =\xf\circ \nf 
= (x_{n_k})_{k \in \N}$ for the corresponding
subsequence of~$\xf$ induced by~$\nf$. 

$\MM(X)$ denotes the set of Borel 
probability measures on~$X$, equipped with the
compact and metrizable topology of weak convergence.
For the special case~$X=[0,1]$ we simply write~$\PP=\MM([0,1])$.   
Sometimes we write $\mu(i)$ for~$\mu(\{i\})$ ($i \in X$, $\mu \in \MM(X)$).  
Let, as usual, $\delta_x \in \MM(X)$ denote the point measure concentrated 
in $x \in X$, i.e.\ $\delta_x(B)=1$ for $x \in B$, $\delta_x(B)=0$ 
for $x \notin B$, $B \tm X$ Borel. 

In order to describe the distribution
behavior of sequences we introduce the discrete measures
\[
\mu_{\xf,N} = \frac 1N \sum_{n=1}^N \delta_{x_n}
\]
and define 
\[
M(\xf):= A((\mu_{\xf,N})_{N \in \N}) \tm \MM(X),
\]
the set of so-called limit measures of~$\xf$.
(\here This set, viewed as set of distribution functions, is called 
$G(\xf)$ in \cite{SP}.\there)
$\xf \sim \lambda$ means that $M(\xf) = \{\lambda\}$, 
i.e.\ the sequence~$\xf$ is {\em uniformly distributed} 
with respect to the measure $\lambda \in \MM(X)$.

The set~$S_0$ of strictly increasing sequences $\nf=\nk$ of positive integers
carries (via $\nf\mapsto \sum_k {2^{-n_k}}$)
a natural measure theoretic structure 
as well as a metric and topological one.
Thus for any given sequence $\xf = \xn$ on some compact metric space~$X$
it makes sense so say that a {\em typical} subsequence 
$\xf \nf = (x_{n_k})_{k \in \N}$ has a certain
property if the set of exceptional $\nf \in S_0$
is small in the sense that either it has
measure zero (measure theoretic), or small Hausdorff dimension 
(depending on the metric), or that it is meager (of first Baire category).
In this paper we will mainly focus on the last, i.e.\ the topological
point of view. Thus we will say that a typical or generic element 
has some property~$P$ if the set of elements with property~$P$  
is residual, i.e., if the set 
of exceptions is meager.
In \cite{GSW} we have investigated the situation with respect to the
set~$M(\xf)$ of limit measures of a sequence,
which is a natural object to describe the distribution
behavior of~$\xf$.

In the measure theoretic
context the typical distribution of a subsequence is the same as for the
original one, i.e.\ $M(\xf \nf)=M(\xf)$. In particular, if $\xf$ is
uniformly distributed w.r.t.\ some measure~$\lambda$, then the same
holds for almost all subsequences (cf.\ also \cite{T} and \cite{LT}).
In the topological context
the situation is quite different, namely: 
Provided all $x \in X$ are accumulation
points of~$\xf$ then $M(\xf \nf)=\MM(X)$ for a generic $\nf \in S_0$.

Note the analogy to the following facts.
Consider the product space $X^{\N}$ of all sequences on~$X$, equipped
with the product measure $\lambda^{\N}$ induced by some fixed 
probability measure $\lambda$ on~$X$. As a consequence of the strong law
of large numbers, $\lambda^{\N}$-almost all $\xf$ are $\lambda$-uniformly
distributed, i.e.~$M(\xf)=\{\lambda\}$, while $M(\xf) = \MM(X)$ for
$\xf \in X^{\N}$ generic in the Baire sense.  Modifying a concept from \cite{M},
sequences $\xf$ with $M(\xf)=\MM(X)$  have been called 
{\em maldistributed} in \cite{Wi}.

Thus the situation is, roughly spoken, as follows: Almost all
(sub)sequences are regularly (uniformly) distributed, but generic
sequences are irregularly distributed  (maldistributed).
The topic of Section 3 is a refined analysis of this topological
maldistribution phenomenon.

\subsection{Kuratowski-Ulam's theorem and $n\alpha$-sequences}\label{SSkurulam}

The theorem of Ku\-ra\-towski-Ulam is the topological counterpart to the
measure theoretic Fubini theorem on product spaces. 
Recall that a Polish space is a complete separable metric space.

\begin{Prop}\label{Pkurulam}{\bf (Kuratowski-Ulam)}
Let $A,B$ be Polish spaces and let $M \tm A \times B$ be a Borel set. 
Furthermore let, for each $a \in A$, $_aM = \{b \in B:\ (a,b) \in M\}$ and
for each $b \in B$, $M_b = \{a \in A:\ (a,b) \in M\}$.
Then the following statements are equivalent:
\begin{enumerate}
\item $M$ is meager in $A \times B$.
\item The set of all $b \in B$ such that $M_b$ is not meager in~$A$
  is meager in~$B$.
\item The set of all $a \in A$ such that $_aM$ is not meager in~$B$
  is meager in~$A$.
\end{enumerate}
\end{Prop}
For proofs and much more background we refer to \cite{O}.

For our context, think about the spaces $A=S_0$ and $B=X=\R/\Z$ (unit 
circle, one dimensional torus).
For each point $(\nf,\alpha) \in S_0 \times X$ we are interested
in $M(\nf\alpha)$ where $\nf \alpha = (n_k\alpha)_{k \in \N}$.
The sequence $\nf\alpha$ is uniformly distributed w.r.t.\ 
Lebesgue (Haar) measure~$\lambda$, hence 
dense in~$X$ for every irrational~$\alpha$.
(Of course $\alpha \in X = \R/\Z$ is called irrational if it is
a remainder class consisting of irrational numbers.)
Theorem 1.3 in \cite{GSW} says that the typical subsequence of a dense
sequence is maldistributed, hence for each irrational $\alpha$
the equality $M(\nf\alpha) = \MM(X)$ holds for a generic~$\nf$.
Since rationals, forming a countable set, are of first category, this shows
that the third condition in Kuratowski-Ulam's theorem is satisfied.
This yields that also the other two conditions hold. The first one
translates to the statement that maldistribution holds for
a generic $(\nf,\alpha) \in S_0 \times X$. The second one, finally,
reads as follows: The set $R \tm S_0$ of all $\nf$
such that the sequence $\nf \alpha$ is maldistributed
for a generic $\alpha$ is residual.
For a fixed $\nf=\nk$ it might be much more difficult to decide
whether it is in~$R$. Section \ref{Smg} will be devoted to this topic,
in particular for sequences satisfying growth conditions.

\subsection{Contents of the paper}\label{SScontents}

The theorem of Kuratowski-Ulam motivates two types of questions. 
They correspond to the main sections of this paper, which can be
read independently of each other.

\setbox0\hbox{Question 1:} 
\begin{list}{-}{\setlength{\leftmargin}{\wd0}%
		\addtolength\leftmargin\labelsep
		\setlength{\labelwidth}{\leftmargin}%
}
\item[Question 1:] Given~$\nf$, can we make assertions on~$M(\nf\alpha)$ for
  generic~$\alpha$? (Section~\ref{Smg})
\item[Question 2:] Given $\alpha$ (or more generally $\xf$ with certain known
  distribution properties), can we make assertions on~$M(\nf\alpha)$ for
  generic~$\nf$? (Section~\ref{Srw} treats a refinement of this question.)
\end{list}

Concerning Question 1, 
it is clear that for sequences $\nf$ with positive lower density
the distribution of~$\nf\alpha$ cannot be arbitrarily irregular.
(However, see \ref{refcomment}.)

This indicates that very strong irregularity results 
(stronger for instance than Theorem \ref{refcomment}) 
can be expected only if the sequence
$\nf$ grows fast enough. A positive result into this
direction is Theorem 1.1 from \cite{salat}: 
If $n_{k+1} =a_k n_k$ with $a_k \in \{2,3,\ldots\}$ for all $k$
then $\nf\alpha$ is not uniformly distributed for generic~$\alpha$.
Our Theorem \ref{salat2} tells us that the same conclusion holds
under the weaker assumption
$\liminf_{k \to \infty} \frac{n_{k+1}}{n_k} > 1$.
Under the stronger growth condition 
$\lim_{k \to \infty} \frac{n_{k+1}}{n_k} = \infty$ 
we can even obtain maldistribution for typical $\alpha$ (Theorem \ref{salat3}).

For arbitrary $\nf$ the situation is not clear. 
We illustrate this by contrasting the cases $n_k = 2^k$ and
$n_k = 2^k+1$ (Theorem \ref{examples}).

Thus the following
very general problem might be an initial point for future research.

{\bf Problem 1:} For which sequences $\nf$ is there a set $M$ of measures
such that $M(\nf \alpha)=M$ for generic~$\alpha$?

The rest of the paper (Section \ref{Srw}) is motivated by Question~2. 
To understand our approach, recall first that Theorem 
1.3 in \cite{GSW} gives a complete answer to the question as stated above: 
Given $\alpha$, a generic subsequence takes as limit measures all 
Borel measures. To get deeper insights we look at appropriate
closed subspaces $S$ of the Baire space $S_0$ of all~$\nf$.
Varying the subspace~$S$ one tries to get different sets $M(S)$ of measures
such that 
$$\hbox to \hsize{\rlap{$(*)$}
\hfill $M(\xf \nf) = M(S)$ for generic $\nf \in S$. \hfill} 
$$
This indeed
works for all $S$ from a certain class of subspaces, each of them
induced by a given interval partition $(I_j)_{j \in \N}$ of $\N$
and a sequence $m_1,m_2,\ldots \in \N$ by the requirement that each $I_j$ 
contains exactly $m_j$ elements from $\nf\in S$.

To save notation at this place we refer to Section \ref{Srw} for
more precise statements. (Note the analogy to stochastic processes
as Markov chains where the probability measure on the space of 
sequences is not the usual product measure but may be supported 
on some small, i.e.\ nowhere dense closed subspace.)

The essential property we will use in the proof is that the sequences
$n\alpha$ are not merely uniformly distributed but even well distributed
(cf.\ \cite{KN} or \cite{DT}). Thus Section 3 will be presented
in this more general context.

Our results (Theorems \ref{Trw1} and \ref{Trw2}) are just first
examples for a topic which might deserve further investigations
in future research. To make such projects more concrete we pose
the following problems:

{\bf Problem 2:} Our results only depend on the well distribution
property but do not make further use of the arithmetic structure of
$n\alpha$-sequences. Thus it seems desirable to find interesting classes 
of subspaces $S$ allowing results of the above type with more
number theoretic impact.

{\bf Problem 3:} Sets $M(S)$ as in $(*)$  cannot exist for arbitrary
closed subspaces $S \tm S_0$. (Every disjoint union $S=S_1 \cup S_2$
with $M(S_1) \neq M(S_2)$ works as a counterexample.) Is it possible
to characterize those $S$ for which there is a set M(S) such that $(*)$ holds?

\section{Sparse subsequences of~$(n\alpha)$}\label{Smg}

In this section $\lambda $ denotes the Lebesgue (Haar) measure on~$\R/\Z$. 

\subsection{Statement of the main  results of this section}\label{SSmg.results}

In this section we consider the distribution behavior of sparse
subsequences of $(n\alpha)_{n\in \N}$. 
In \cite{salat},  (essentially) the following has been proved: 

\begin{Prop} \label{salat1} {\bf ({\v S}al{\'at})}
Let $\nf=(n_0,n_1,\ldots)$ be a sequence of natural numbers satisfying 
$n_{k+1} \ge 2n_k$ for all~$k$.  Then the set 
$${\mathscr U}:= \{\,\alpha\in \R/\Z: \nf \alpha \mbox{ is uniformly 
distributed w.r.t.~$\lambda$}\,\}$$
is meager. 
\end{Prop}

We will improve this result by weakening the growth condition on the
sequence~$\nf$.

 \begin{Def}
For any sequence $\xf = (x_n)$ and any interval~$I$, 
we define $\bar\mu_\xf(I)$
by 
$$ 
\bar\mu_\xf(I): = \sup\{\,\mu(I): \mu\in M(\xf)\,\}.$$
\end{Def}

\begin{Rem}\label{muquer}
Note that $\bar\mu_\xf(I) \ge \limsup_{n \to \infty} \mu_{\xf,n}(I)$ 
while equality does not hold in general: Take $x_n = \frac 1n$ and
$I=\{0\}$, then~$M(\xf)=\{\delta_0\}$, $\delta_0(I)=1$ but
$\mu_{\xf,n}(I) = 0$ for all~$n$.   \here It is also easy to see that 
$\bar \mu$ is in general not additive, hence not a measure.\there
\end{Rem}

 \begin{Theo}\label{salat2}
Let $\nf=(n_0,n_1,\ldots)$ be a sequence of natural numbers, and assume 
$q:=\liminf_k (n_{k+1}/n_k) > 1$.  Then the set 
$${\mathscr U}:= \{\,\alpha\in \R/\Z: \nf \alpha \mbox{ is 
  uniformly distributed w.r.t.~$\lambda$}\,\}$$
is meager. 

Moreover: There is a number $Q>0$ such that for all intervals $I$ 
\here of length $<\frac1q$ \there  the set 
$$
\{\,\alpha: \bar\mu_{\nf \alpha}(I) > \frac{Q}{-\log \lambda(I)}\,\}
$$
is residual. 

Equivalently, the set 
$\{\alpha:  \forall I\, \bar\mu_{\nf \alpha}(I) > \frac{Q}{-\log \lambda(I)}\}$
is residual \here (where the quantifier $\forall I$ refers to all intervals of length $<1/q$)\there.
\end{Theo}

\begin{Rem}\label{mg.rem}
\begin{enumerate}
\item The sentence ``Equivalently \dots''\ follows from the previous
 sentence because it is enough to prove this 
 for intervals with rational end points.   
\item Note that for short intervals $I$ we 
have $\frac{Q}{-\log \lambda(I)} \gg \lambda(I)$. 
\here\there
% \item Choosing $Q$ small enough, the inequality 
% $\bar \mu_{\nf \alpha}(I) > \frac{Q}{-\log \lambda(I)}$ 
% will be trivially true for large
% intervals~$I$, say for all $I$ with $\lambda(I)\ge \frac 1 q $. 
% So it is enough to consider only intervals $I$ with 
% $\lambda(I) < \frac 1 q $. 
\item Results from \cite{AHK} or \cite{B} show that the growth condition 
in Theorem \ref{salat2} cannot be weakened. Boshernitzan for instance shows
that for every sequence of integers $m_k$ with 
$\lim_{k \to \infty} % m_k^{\frac 1k}
 \sqrt[k]{m_k\vphantom{I}} = 1$ there are $n_k \ge m_k$
such that $(n_k\alpha)_{k \in \N}$ is uniformly distributed mod 1 for
all irrational~$\alpha$.
\end{enumerate}
\end{Rem}

{\sl \noindent
Idea of the proof of Theorem~\ref{salat2}: Fix a short interval~$I$. 
Let $c$ be large with respect to~$q$ and $I$  (see below for details).
  If we consider
only every $c$-th term in the sequence~$\nf$, i.e., the sequence
$\nf' = (n'_k)_{k\in \N}$ with $n_k' = n_{c k}$, then the $n'_k$ will 
increase so fast that 
\begin{equation}\label{sketch}
\hbox{\rm $R:= \{\alpha: $ $\{n'_k\alpha, \ldots,
  n'_{2k-1}\alpha\}\subseteq I$ for infinitely many $k$ $\}$
is residual}.
\end{equation}
So for $\alpha\in R$, $\bar\mu_{\nf'\alpha}(I)\ge \frac12$, and 
 $\bar\mu_{\nf\alpha}(I) \ge  \frac{1}{2c}$. 
\newline
Upon closer inspection we see that $c\approx \frac{1}{-\log_q {\lambda(I)}}$ 
is sufficient for (\ref{sketch}).
\newline
Similar methods will be used in the proof of Theorem~\ref{salat3}.

}

 \begin{Theo}\label{salat3}
Let $\nf=(n_0,n_1,\ldots)$ be a sequence of natural numbers, and assume 
$\lim_k (n_{k+1}/n_k) = \infty$.  Then the set 
$$ \{\,\alpha\in \R/\Z: \nf \alpha \mbox{ is maldistributed}\,\}$$
is residual. 
\end{Theo}

Weaker versions of irregular distribution also occur for certain
classes of slowly increasing $n_1 < n_2 < \ldots$. The following
two theorems elaborate on remarks of the referee of a previous
version of this paper, for which we are grateful:

\begin{Theorem}\label{refcomment}
Let $X$ be the set of all increasing $\nf = (n_k)_{k \in \N}$ of integers 
with $n_{k+1}-n_k\in \{1,2\}$.  
Then the set $\{(\alpha, \nf)\in \R \times X:
\mbox{$\nf \alpha$ is not u.d.}\}$  
is  residual in $\R\times X$. 
\end{Theorem}

More refined investigations in this spirit will be the content
of Section \ref{Srw}.
 
The last result of this section indicates that
for individual sequences of $\nf$ the situation can be very diverse
and hence complicated:

\begin{Theorem}\label{examples}
\here For $\nf = (2^k)_{k \in \N}$ there is a set $M \tm \MM([0,1])$ such that 
$M = M(\nf \alpha)$ for generic $\alpha \in [0,1]$. 
$M$ contains exactly all measures which are invariant under
$x \mapsto 2x$. In contrast,
for $\nf' = (n_k')_{k \in \N}$ with $n_k' = 2^k+1$, there is no $M' \tm \MM([0,1])$
such that $M'=M(\nf'\alpha)$ for generic $\alpha \in [0,1]$.\there
\end{Theorem}

For related constructions yielding results in terms of Hausdorff dimension
we refer to \cite{P}.

\subsection{Notation}\label{SSmg.notation}

For notational convenience we sometimes identify $\alpha + \Z \in \R/\Z$
with the unique representative $\alpha \in [0,1) \tm \R$.
Very often we are in the situation
that an intersection $I \cap B$ of an interval $I$ with a Borel set $B$
is residual in~$I$. Note that this can be interpreted as a generalized
implication of the type:
\begin{quote} Except for a 
 meager set, $x \in I$ implies $x \in B$.
\end{quote}
Therefore we introduce the following notation.

\begin{Def} For an open interval $I$ and a Borel set $B$ we write 
$$ I\forces B$$
as an abbreviation for ``$B\cap I$ is residual in~$I$'' or
equivalently, ``$I\setminus B$ is meager''.   We read this also
as
``the typical element of~$I$ is in~$B$''.
\end{Def}

The following fact is a folklore consequence of Baire's theorem: 

\begin{Fact}\label{f-}
Let $I$ be an open interval. 
\begin{enumerate}
\item  If~$B_n$ is a Borel set for every
$n\in \{0,1,2,\ldots\}$, and 
 $ I \cap \bigcup_n B_n $ is residual  in~$I$, then there is some 
open nonempty $J \subseteq I$ and some $n$ such that $B_n$ is residual in~$J$, 
or abbreviated: 
$$ I  \forces \bigcup _{ n\in \N}  B_n \ \Rightarrow  \
\exists J \subseteq I\, \exists n\in \N: \ J \forces B_n$$
\item  If~$B_n$ is a Borel set for every
$n\in \{0,1,2,\ldots\}$, then    
 $ I \cap \bigcap_n B_n $ is residual in~$I$ iff each $I \cap B_n$ 
is residual in~$B_n$:   
$$ I  \forces \bigcap _{ n\in \N}  B_n \ \Leftrightarrow  \
\forall n\in \N: \ I \forces B_n$$
\item If~$B$ is a Borel set then $B\cap I$ is {\em not}
 residual in~$I$ iff there is
 some open interval $J \subseteq I$ such that $B$ is meager in~$J$: 
$$ I\nVdash B  \ \Leftrightarrow \ \exists J \subseteq I: J \forces (-B).$$ 
(Here we write $-B$ for the complement of~$B$.) 
\end{enumerate}
\end{Fact}

\newcommand{\B}{{\bf B}}
\newcommand{\symdiff}{\mathbin\Delta}
\begin{proof}
Let $\B$ be the family of all sets with the Baire property, 
i.e., all sets which
can be written as $A \symdiff M$, where $A$ is an open set,
$M$ is meager, and $\Delta$ denotes the symmetric difference of two sets. 
Then clearly 
\begin{itemize}
\item $\B$ contains all open sets
\item $\B$ is closed under countable unions
\item $\B$ contains all closed sets, as each closed set $C$
can be written as $A \symdiff M$, where $A$ is the open kernel of $C$ and 
$M = \partial C = C  \setminus A$ is nowhere dense.
\item $\B$ is closed under complements: If $X = A  \symdiff M$, then 
$(-X) = (-A) \symdiff M$; write $-A$ as $A' \symdiff M'$ with $A'$
open, $M'$ meager, then $(-X) = (A' \symdiff M') \symdiff M =   
A' \symdiff( M' \symdiff M )$, where  $M' \symdiff M$ is meager. 
\end{itemize}
Hence $\B$ contains all Borel sets.  

To prove (1), write each $B_n$ as $A_n\Delta M_n$ with $A_n$ open 
and $M_n$ meager. Not all $A_n$ can be empty (otherwise the set $\bigcup_n B_n
= \bigcup_n M_n$ would be meager); let $J$ be an interval contained in
any nonempty~$A_n$.

(2) is easy. 

To prove (3), assume that $B\cap I$ is not residual, and write 
 $I\setminus (B\cap I)$ as $ A \symdiff M$ for some open $A$ and meager~$M$;
as $I\setminus (B\cap I)$ is not meager, $A$ is not empty; 
let $J$ be any nonempty open interval with $J\subseteq A$. 
\end{proof}

\begin{Def}
We say that a family $(f_1, f_2,\ldots)$ of functions $f_i: [0,1] \to [0,1)$
is $\varepsilon $-mixing if: whenever $J_1,J_2,\ldots$ are intervals 
of length~$\varepsilon$, then for all $k\in \N$:  $$ \bigcap_{n=1}^k
 f^{-1}_n(J_n) \ \mbox{ contains an inner point.} $$
More generally we say that  $(f_1, f_2,\ldots)$ 
is ${\ve}$-mixing in~$\delta$ if: for all 
sequences $J_1,J_2,\ldots$ of intervals of
length $\varepsilon $,  and all $k\in N$, and all intervals~$J'$ 
of length $ {\delta} $:  
\[ 
  J' \cap \bigcap_{n=1}^k f^{-1}_n(J_n) \ \mbox{ contains an inner point.} 
\]
\end{Def}

\begin{Rem}
Although we work in~$I=[0,1)$, we identify elements in~$I$ with their
equivalence classes in~$\R/\Z$, so an (open) interval can either be of
the form $(a,b)$ or of the form $[0,a) \cup (b,1)$ (for $0\le a<b\le
1$). However, since we are mainly concerned with very short intervals,
it is no loss of generality to only consider intervals of the first
form.
\end{Rem}

\subsection{Proof of Theorem \ref{salat2}}\label{SSproof.salat2}

\begin{Lemma}\label{17} 
Let $f_k:[0,1)\to[0,1)$ be the function mapping~$ \alpha  $ to~${n}_k \alpha  $
modulo~$1$, where $\nf= \here(n_k)_{k=1}^\infty\there$ is a sequence of 
natural numbers satisfying
\begin{enumerate}
\item $ {n}_{k+1} > \frac{2}{ \varepsilon } {n}_k$ for all~$k\ge 1$.
\item $ {n}_1 >\here \frac{2} {{\delta} }\there$. 
\end{enumerate}
Then $(f_1, f_2,\ldots\, )$ is $ {\ve}$-mixing in~$\delta$. 
\end{Lemma}

\begin{proof}
\here For notational convenience let $n_0$ be a real number satisfying $
\frac\varepsilon\delta < n_0 < \frac\varepsilon2 n_1$. So we have 
$$ \delta >\frac \varepsilon{n_0} \mbox{ \quad and \quad} \frac\varepsilon{n_{k-1}} > 
\frac 2{n_k} \mbox { for $k=1,2,\ldots$}\there$$
Let~$J_1, J_2,\ldots$ be intervals of length 
$ \varepsilon $,  $ J' $ an interval of length $ {\delta} $.

 We will show (by induction on~$k\here =0,1,2,\ldots\there$) that
each set  $$ J' \cap   \bigcap_{i=1}^k f^{-1}_i(J_i)$$ 
contains in fact an interval $I_k$ of length $ \varepsilon/{n}_k$. 
This is clear for~$k=0$, as the length of~$J'$ is $ {\delta} >
\varepsilon/{n}_0$. 

Consider~$k>0$. Note that $f_k^{-1}(J_k) = 
 \bigcup_{j=1}^{{n}_k} \{\alpha\in [0,1) : {n}_k \alpha  - j \in J_k \}$ 
is a union
 of~${n}_k$ many disjoint intervals, each of length $ \varepsilon /{n}_k$. 

By inductive assumption, the set $ J' \cap   \bigcap_{i=1}^{k-1}
f^{-1}_i(J_i)$ contains an interval $I_{k-1} $ of
 length $ \varepsilon/{n}_{k-1}  $:
$$ I_{k-1} \subseteq 
 J' \cap   \bigcap_{i=1}^{k-1}
f^{-1}_i(J_i), \qquad  \lambda(I_{k-1}) = \frac{\e}{n_{k-1}}$$

 Since  $ \varepsilon / {n}_{k-1} > 2/{n}_k$, 
we can find a natural number~$j < {n}_k$
such that the interval  $$ [\frac{j}{{n}_k}, \frac{ j+1}{{n}_k} ] =
 \{\, \alpha : {n}_k  \alpha  - j \in
[0,1]\,\}$$
is contained in~$I_{k-1} $. 
Hence the set  
\[I_k:= \{\,\alpha : {n}_k \alpha  - j \in J_k \,\}, 
\]
an interval of length $ \varepsilon  / {n}_k$, is also contained
 in~$I_{k-1}$. 
\end{proof}

%\end{mixing}

\begin{proof}[Proof of Theorem~\ref{salat2}]  
\here We will prove the ``moreover'' statement in Theorem~\ref{salat2}.\there

Choose $Q>0$ so small that \ $(\frac{1}{4Q}-1)-\log 2 > 1$. 

Without loss of generality we may assume $\forall k: \frac{n_{k+1}}{n_k} > q 
\here > 1 \there$.

Let $\ve:= \lambda(I)<\frac{1}{q}$,\here\there
so $(-\log \ve) > 1$.   (In this proof, $\log$ denotes the logarithm with 
base~$q$.)

 So we have
                $(\frac{1}{4Q}-1)\cdot(-\log \ve)-\log 2 > 1$,
hence the interval 
$$ (\log 2 -\log \ve, -\frac{1}{4Q} \log \ve) $$
has length~$>1$.   Let $c$ be an integer in this interval. 
Thus,
\begin{itemize}
\item $q^c > \frac 2 \ve$  
\item $\displaystyle \frac 1 {2c} > \frac{2Q}{- \log \ve}$
\end{itemize}

Now assume that the theorem is false.   Since the set 
$\{\alpha: \bar\mu_{\nf \alpha}(I) > \frac{Q}{-\log\ve}\}$
is a Borel set and not residual, by \ref{f-}(3) we know that its
complement will be residual
in~$\here J_0\there$, for some open interval~$J_0$:  
$$J_0 \forces 
\biggl\{\,\alpha: \bar\mu_{\nf \alpha}(I) \le  \frac{Q}{-\log\ve}\,\biggr\}$$

Now, by Remark~\ref{muquer}, the set 
$ \bigl\{\,\alpha: \bar\mu_{\nf \alpha}(I) \le  \frac{Q}{-\log\ve}\,\bigr\}$
is contained in the set 
$$
\biggl\{\,\alpha: \exists m \,\forall N\ge m: \mu_{\nf \alpha,N}(I)  < 
\frac{2Q}{-\log\ve} \,\biggr\}.
$$ 

We will write $\anzahl{N}{\alpha}$ for the set
 $\{ j<N : {n}_j \alpha \in I \}$. So $\mu_{\nf\alpha,N}(I) = 
\frac{\#\anzahl{N}{\alpha}}{N}$ and we have 
$$ J_0 \forces 
\bigcup_m \bigcap_{N\ge m}
	\biggl\{\, \alpha: \frac {\#\anzahl{N}{\alpha}}{N} <
 \frac {2Q}{-\log\ve} \,\biggr\}.$$
By \ref{f-}(1), 
 we can find an open interval
 $\here J_1\subseteq J_0\there$ and a $k^*$ such that 
$$ 
 J_1 \forces 
 \bigcap_{N\ge k^*}
	\biggl\{\, \alpha: \frac {\#\anzahl{N}{\alpha}}{N}  < 
	\frac {2Q}{-\log\ve} \,\biggr\}.
$$
In other words: for all $N\ge k^*$: 
\begin{equation}\label{res}
 J_1 \forces 
 	\biggl\{\, \alpha: \frac {\#\anzahl{N}{\alpha}}{N}  < 
	\frac {2Q}{-\log\ve} \,\biggr\}.
\end{equation}

Let $\delta:= \lambda(J_1)$.  
We have $n_{k^* c} > \frac{2}{\delta}$ (otherwise
we just increase $k^*$).

Now we consider the functions $f_{k^* c}, f_{(k^*+1)c},  f_{(k^*+2)c}, 
\ldots, f_{(2 k^*-1) c }$.   Since 
$$ 
\frac{n_{(k^*+i+1)c}}{n_{(k^*+i)c}}\ge q^c 
 > \frac 2 \ve,
$$
and $n_{k^* c } > \frac{\ve}{\delta}$, these functions are
$\ve$-mixing in $\delta $ (Lemma \ref{17}). 

So there is an open interval $\here J_2 \subseteq J_1\there$ such that for all 
$\alpha\in J_2$, and all $i\in \{0,\ldots, k^*\}$:
$$
\alpha\in f^{-1}_{n_{(k^*+i)c}}(I) \qquad \mbox{i.e., \ \ }
n_{(k^*+i)c} \alpha \in I.$$

Thus $\forall \alpha\in J_2$:
$$
 \#  \anzahl{2k^*c}{\alpha}  =  \#  \{i < 2k^*c: n_i \alpha\in I\} \ge \# 
\{\, k^*c, (k^*+1)c,\ldots, (2k^*-1)c\,\} = k^*.
$$ 

Hence  for $\alpha \in J_2$:
\begin{equation}
\label{empty}
 \frac{\#\anzahl{2k^* c}{\alpha}}{2 k^* c} > \frac{1}{2c}.
\end{equation}

However, $\frac{1}{2c} >  \frac{2Q}{-\log\ve}$ and $J_2 \subseteq J_1$,
 so we get from (\ref{res}) for $N:= 2k^*c$:
\begin{equation}\label{res2}
J_2 \forces 
	\biggl\{\, \alpha: \frac{\#\anzahl{2k^*c}{\alpha}}{2k^*c}
                     \le \frac{1}{2c} \,\biggr\}.
\end{equation}
Now consider the set 
	$\bigl\{\, \alpha: \frac{\#\anzahl{2k^*c}{\alpha}}{2k^*nc}
                     \le \frac{1}{2c}\, \bigr\}\cap J_2$.
By (\ref{empty}), this set is empty, but by (\ref{res2}) it is residual
in~$J_2$; this is a contradiction.

\end{proof}

\subsection{Proof of Theorem \ref{salat3}} \label{SSproof.salat3}

\begin{Fact}\label{closed1}
For any sequence $\xf = \xn$, the set $M(\xf)$ is closed. 
\end{Fact}

\begin{Fact}\label{closed2}
For any sequence $\nf = \nk$, the set 
$$ M:= \{ \mu\in \PP:
\mbox{ for typical $\alpha$, $\mu\in M(\nf \alpha)$} \,\}$$
is closed in $ \PP =  \MM([0,1])$.
\end{Fact}
\begin{proof} We show that $M$ is closed under limits of sequences. 
So let $\mu_n \to \mu$, with all $\mu_n\in M$.  Let 
$$ A_n:= \{\alpha: \mu_n\in M(\nf \alpha)\}$$
Now $\mu_n\in M$ just means that $A_n$ is residual; so
$A:= \bigcap_n A_n$ is also residual, and by \ref{closed1} 
we have $\mu\in M(\nf \alpha)$ for all $\alpha\in A$.  
\end{proof}

\begin{Def}
For any list $\vec e = (e_0,\ldots, e_{\ell-1})$ of natural numbers, and 
any $\eta>0$  we let 
$$ M_{\vec e, \eta}:= \biggl\{\mu\in \PP:\ \forall i\in \{0,\ldots,\ell-1\}\ 
\bigl|\mu(\halfopen{\frac{i}{\ell}}{\frac{i+1}{\ell}})
              - \frac{e_i}{e} \bigr|<\eta\biggr\}
$$
and $e:= \sum e_i$.

\end{Def}
By \ref{closed2}, the following are equivalent for any~$\nf$: 
\begin{enumerate}
\item[(i)]
     The set $\{\alpha: M(\nf\alpha) = \PP \}$ is residual.
\item[(ii)]
 For each $\vec{e}$ and each~$\eta$, the set 
 $\{\alpha: M(\nf\alpha)\cap M_{\vec e, \eta} \not=\emptyset \}$ is residual.
\item[(iii)]
 For each $\vec{e}$ and each~$\eta$, the set 
 $\{\alpha: \exists^\infty N\,  \mu_{\nf\alpha, N} \in
M_{\vec e, \eta}\}$ is residual.
\end{enumerate}

\begin{proof}[Proof of Theorem~\ref{salat3}] 
Assume that property (iii) above fails. As in the proof of~\ref{salat2}, this
means that we can find a nonempty interval~$I$, a natural number~$N_0$, 
a sequence $\vec e = (e_0,\ldots, e_{\ell-1})$ of natural numbers, 
and a real number $\eta$ such that 
\begin{equation}\label{s3.1}
 I \forces \{\alpha: \forall N\ge N_0: \mu_{\nf\alpha,N}\notin 
 M_{\vec e, \eta} \}.
\end{equation}
Clearly we may assume~$N_0>\frac{1}{\eta}$, 
that $e:= \sum e_i$ divides~$N_0$, and that 
\begin{equation}\label{s3.1a}
\frac{n_{k+1}}{n_k} > 2 \ell, \quad n_{N_0} > \frac{1}{\lambda(I)}.
\end{equation} 

   Choose a sequence $(I_j: j=1,\ldots, N_0^2)$
of intervals such that for each $i\in \{0,\ldots, \ell-1\}$ the 
set 
$$ \{ j \in \{1,\ldots, N_0^2\}: \ I_j = \halfopen{\frac{i}{\ell}}
 {\frac{i+1}{\ell}} \}$$
has cardinality $\frac{e_i}{e}N_0^2$. So each $I_j$ 
has length $\frac 1{\ell}$.

Let $f_j (x) = n_j x$ for $j\in \{N_0+1,\ldots, N_0^2\}$.  
By (\ref{s3.1a}) and Lemma \ref{17} these
functions are $\frac{1}{\ell}$-mixing in~$\lambda(I)$, 
so we can find an interval 
$$ J \subseteq I \cap \bigcap_{j=N_0+1}^{N_0^2} f^{-1}_j(I_j).$$
We now claim that 
\begin{equation}\label{s3.2}
\forall \alpha\in J: \mu_{\nf\alpha, N_0^2} \in M_{\vec e,\eta}, 
\end{equation}
which clearly contradicts (\ref{s3.1}). 

Indeed, let $\alpha\in J$.  Then for any $j\in \{N_0+1,\ldots, N_0^2\}$ 
we have $f_{j}(\alpha) \in I_j$, hence 
(writing $O(1)$ for a quantity that lies between $-1$ and $1$) we get
$$
\mu_{\nf\alpha, N_0^2} (\halfopen{\frac{i}{\ell}}{\frac{i+1}{\ell}}) = 
\frac{1}{N_0^2}\bigl( \frac{e_i}{e}N_0^2 + O(1) N_0\bigr) 
= \frac{e_i}{e} + \frac{O(1)}{N_0} , $$ 
so $ \mu_{\nf\alpha, N_0^2} \in M_{\vec e, \eta}$.  
\end{proof}

\subsection{Proof of Theorem \ref{refcomment}} \label{SSproof.refcomment}

The set in question has the Baire property, and we will show that 
each vertical section is residual.

So fix~$\alpha\in(0,1)$. We can find a number $\varepsilon>0$ such that 
the intervals $(0,\varepsilon)$ and $(\alpha, \alpha+\varepsilon)$
(computed modulo 1) are disjoint.

We claim
\begin{quote}
Whenever 
 $(n_1,\ldots, n_k)$ is a finite sequence in which $n_{j+1}-n_j\in \{1,2\}$ 
 holds for all~$j<k$,
\\
there is an infinite extension $(n_1,\ldots, n_k, n_{k+1},
n_{k+2},\ldots)\in X$ such that 
$$ \forall j>k:  \ n_j\alpha\notin (0,\varepsilon).$$
\end{quote}

This claim implies that for each $k$ the closed set 
$$ \bigcap_{j>k} \{ \vec n:   
\frac{|\{i<j: n_i\alpha\in (0,\varepsilon)\}|}{j} \ge \frac{\varepsilon}{2}\}$$
is nowhere dense, so the set $\{\vec n: \mbox{$\vec n \alpha$ is u.d.}\}$
is meager.

Proof of the claim:  We can construct the numbers $n_{k+1}, n_{k+2},\ldots
$ by induction.  
Given~$n_{j}$, we either have
$(n_j+1)\alpha \notin (0,\varepsilon)$ --- in that case 
we may choose $n_{j+1}:= n_j+1$.   Or we have 
$(n_j+1)\alpha\in (0,\varepsilon) $ --- in that case  we have 
$(n_j+2)\alpha\in (\alpha,\alpha+\varepsilon)$, hence 
$(n_j+2)\alpha\notin (0,\varepsilon)$, so we may 
 choose $n_{j+1}:= n_j+2$.

\subsection{Proof of Theorem \ref{examples}} \label{SSproof.examples}

For Theorem \ref{examples} it suffices to prove the four statements
of the following lemma. Recall that we focus on the sequences 
$\nf = (n_k)_{k \in \N}$ with $n_k = 2^k$ and $\nf' = (n_k')_{k \in \N}$
with~$n_k'=2^k+1$.

\begin{Lemma}\label{Lexamples}
\begin{enumerate}
\item Let $X$ be any compact metric space, $T: X \to X$ continuous,
$x \in X$, $\xf = (T^nx)_{n \in \N}$ and $\mu \in M(\xf)$. 
Then $\mu$ is $T$-invariant.
\item Let $\mu \in \MM(X)$ be $T$-invariant for $T: X \to X$, $x \mapsto 2x$, 
on~$X=\R/\Z$. Then $\mu \in M(\nf \alpha)$ for generic $\alpha \in X$.
\item If $\alpha \in (0,\frac 1{16})$ then $\mu \in M(\nf' \alpha)$
implies $\mu(I) \le \frac 56$ for $I:=(\frac 12, \frac 34)$.
\item For generic $\alpha \in (\frac 12, \frac 34)$ there is a 
$\mu \in M(\nf'\alpha)$ with $\mu(I) = 1$ for $I:=(\frac 12,\frac 34)$.
\end{enumerate}
\end{Lemma}
\begin{proof}[Proof of Theorem~\ref{examples}]
Assume that Lemma \ref{Lexamples} holds.
Let $M$ denote the set of all $T$-invariant measures 
$\mu \in \MM(\R/\Z)$ for $T: x \mapsto 2x$. Then the first statement of the
lemma tells us that $M(\nf \alpha) \tm M$ for all $\alpha \in X$. Conversely,
the second statement guarantees that for each $\mu \in M$ the set
$R_{\mu} = \{\alpha:\, \mu \in M(\nf\alpha)\}$ is residual. There is
an at most countable set $M_0 = \{\mu_n:\, n \in \N\}$ with~$\overline{M_0}=M$.
Let $R = \bigcap_{n \in \N} R_{\mu_n}$. Then $R$ is residual and
$M_0 \tm M(\nf \alpha)$ for all $\alpha \in R$. Since every set of
the form $M(\xf)$ is closed we have $M = \overline{M_0} \tm M(\nf\alpha)$
for all such~$\alpha$, hence $M(\nf\alpha)=M$ for residual $\alpha \in X$,
establishing the first two sentences in Theorem~\ref{examples},
while the third sentence follows by combining the third and the fourth
statement of the Lemma. 
Thus Theorem \ref{examples} indeed follows from Lemma~\ref{Lexamples}.
\end{proof}

We are now going to prove the four statements of Lemma~\ref{Lexamples}.

\begin{proof}[Proof of  statement (1) of Lemma~\ref{Lexamples}] 
All we have to prove is  $\int f \, d\mu = \int f \circ T \, d\mu$
for any continuous $f: X \to \R$. $\mu \in M(\xf)$ means that 
$\lim_{k \to \infty} \mu_{\xf,n_k} = \mu$ for some $n_1 < n_2 < \ldots \in \N$.
Thus we easily obtain
\[
\begin{array}{rl}
\int f\, d\mu =& \lim_{k \to \infty} \frac 1{n_k} \sum_{j=1}^{n_k}f(T^jx) =
   \lim_{k \to \infty} \frac 1{n_k} \sum_{j=2}^{n_k+1}f(T^jx) = \\
   =& \lim_{k \to \infty} \frac 1{n_k} \sum_{j=1}^{n_k}f \circ T(T^jx) =
   \int f \circ T\, d\mu.
\end{array}
\]
\end{proof}

\begin{proof}[Proof of the  statement (2) of Lemma~\ref{Lexamples}]
Birkhoff's ergodic theorem guarantees that, for each continuous $f: X \to \R$,
the set $R_f$ of all $\alpha \in X$ with 
\[
\lim_{k \to \infty} \frac 1k \sum_{j=1}^k f(2^j \alpha) =
\lim_{k \to \infty} \frac 1k \sum_{j=1}^k f(T^j \alpha) = \int f\, d\mu
\]
has full measure~$\mu(R_f)=1$. Let $f_1,f_2,\ldots$ be any sequence
of continuous $f:X \to \R$, then $\mu(R_{\mu})=1$
for $R_{\mu}: = \bigcap_{n \in \N} R_{f_n}$, in particular 
$R_{\mu} \neq \emptyset$. We may assume that the $f_i$, $i \in \N$,
are $\|{\cdot}\|_{\infty}$-dense in the space of all continuous $f: X \to \R$.
Take any $\alpha_0 \in R_{\mu}$. 
Then $M(\nf\alpha_0) = M((2^k\alpha_0)_{k \in \N}) = \{\mu\}$.
In order to show $\mu \in M(\nf \alpha)$ for generic $\alpha$ assume that
$I \tm X = [0,1)$ is any nonempty open subset of~$X$. Then $I$ contains
an interval $I_0 = (\frac{k-1}{2^{k_0}},\frac k{2^{k_0}})$ for some
$k_0 \in \N$ and $k \in \{1,2,\ldots,k_0\}$. $X$ is a compact metrizable space,
hence so is $\MM(X)$ and there is a sequence 
$U_1 \supset U_2 \supset \ldots$ of open $U_i \tm \MM(X)$ forming
a neighbourhood base for~$\mu$. Note that in $X$ we have
$2^{k_0+k}(\alpha') = 2^k\alpha_0$ for
$\alpha' := \frac{k-1}{2^{k_0}} + \frac {\alpha_0}{2^{k_0}}$.
It follows that for each $j$ there is 
an open neighbourhood $V = (\alpha'-\e_j, \alpha' + \e_j) \tm I_0$, 
$\e_j>0$, such that $\mu_{\nf\alpha,k_j} \in U_j$ for some $k_j \ge j$ 
and all $\alpha \in V$. In particular each
$R_j = \{\alpha:\, \mu_{\nf \alpha,k} \in U_j \ \mbox{for some $k \ge j$}\}$
contains an open dense set. Thus $R = \bigcap_{j \in \N} R_j$ is residual
with $\mu \in M(\nf\alpha)$ for all $\alpha \in R$.
\end{proof}

\begin{proof}[Proof of  statement (3) of Lemma~\ref{Lexamples}]
Assume $0 < \alpha < \frac 1{16}$. Note that for every $J \tm X$ we have
$n_k'\alpha = 2^k\alpha + \alpha \in J$ if and only if 
$n_k\alpha = 2^k\alpha \in J-\alpha$. In order to obtain the desired
estimate we take for $J$ instead of $I=(\frac 12,\frac 34)$ 
the interval $I':= (\frac 12-\frac {\alpha}3, \frac 34 + \frac{\alpha}3)$
and, accordingly $I_{\alpha}:= I'-\alpha = (\frac 12 - \frac 43\alpha,
\frac 34 - \frac 23 \alpha) = I^- \cup I^+$ with 
$I^-: = (\frac 12 - \frac 43\alpha,\frac 12]$ and 
$I^+: = (\frac 12, \frac 34 - \frac 23 \alpha)$. Observe that for
$T:x \mapsto 2x$ each of the sets $T(I^-), T^2(I^-)$ and $T(I^+)$ 
has empty intersection with~$I_{\alpha}$. If $D^-$ and $D^+$ denote the sets
of all $k \in \N$ such that $2^k \alpha \in I^-$ resp.\ $2^k \alpha \in I^+$, 
this shows that the upper densities of $D^-$ and $D^+$ are at most $\frac 13$
resp.\ $\frac 12$. It follows that the set of all $k \in \N$ with
$n_k\alpha = 2^k \alpha \in I_{\alpha}$
or, equivalently, $n_k'\alpha \in I'$
has upper density at most $\frac 12 + \frac 13 = \frac 56$.
Since the interval $I'$ is open and contains the closure of $I$ this yields
that $\mu(I) \le \frac 56$ for every $\mu \in M(\nf'\alpha)$.
\end{proof}

\begin{proof}[Proof of statement (4) of Lemma~\ref{Lexamples}]
Similar arguments as several times before show that
a generic $\alpha$ contains extremely long blocks of 0's in its binary
representation $\alpha = \sum_{j=1}^{\infty} \frac{a_j}{2^j}$.
To be more precise, the set $R$ of all $\alpha$ such that
$(a_j,\ldots,a_{j^2})=(0,\ldots,0)$ for infinitely many $j \in \N$ is residual.
For $\alpha \in (\frac 12, \frac 34) \cap R$ this implies that
the upper density of the set $D = \{k \in \N:\, 
2^{k+1}\alpha = 2^k \alpha + \alpha \in (\frac 12, \frac 34)\}$
is 1, implying that $\mu(I)=1$ for some $\mu \in M(\nf'\alpha)$.
\end{proof}

\section{Baire spaces of subsequences}\label{Srw}

In this section $\lambda $ denotes an arbitrary but fixed Borel probability 
measure on a compact metric space~$X$. 

\subsection{Notation and statement of the main results of this 
  section} \label{SSrw.results}

 We assume $\xf = \xn \approx \lambda$ (see below).
Typical 
examples of this type are $x_n = n\alpha$, $\alpha \in \R/\Z$, or,
more generally,
sequences induced by uniquely ergodic dynamical systems as
ergodic group rotations, i.e.\ $x_n=ng$ where $g$ is a topological
generator of a monothetic compact group, $\lambda$ the Haar measure.
To state our results we need a lot of notation. Therefore the
following list might be for the reader's convenience.

\begin{itemize}
\item We fix a measure $\lambda \in \MM(X)$ (e.g.\ the Lebesgue measure).
\item
$\CC(\lambda)$ denotes the system of all $\lambda$-continuity sets~$C$,
 i.e.\ of those $C \tm X$ with $\lambda(\boundary C) = 0$, where
$\boundary C$ is the topological boundary of~$C$. 
  Similarly $\CC(\lambda,\mu) = \CC(\lambda) \cap \CC(\mu)$ etc. 
\item
$\sigma (\xn) = (x_{n+1})_{n \in \N}$, the shift acting on
arbitrary infinite sequences.
\item
We write $\xf \approx \lambda$ if $\xf \sim \lambda$ in fact is
well distributed. This, by definition, means that 
$\lim_{N \to \infty} \mu_{\sigma^k(\xf),N} = \lambda$ uniformly in $k \in \N$. 
(Since $\MM(X)$ is compact there is a unique uniform structure and
this notion is well defined.)
Clearly $\xf \approx \lambda$ implies $\xf \sim \lambda$
but not conversely. $\xf \approx \lambda$ is equivalent to the condition
that for all $A \in \CC(\lambda)$ the limit 
\[
\lim_{N \to \infty} \frac 1N |\{ k < n \le N+k:\ x_n \in A\}| = \lambda(A)
\]
is uniform in $k \in \N$.
\item
$\If = \Ij$ denotes a  partition of~$\N$ into intervals: 
$I_j = \{n \in \N: a_{j-1} < n \le a_j\}$, $0=a_0 < a_1 < \cdots \in \N$.
The $b_j = a_j-a_{j-1}$ are called block lengths.
\item
For a sequence $\mf = (m_j)_{j \in \N}$ of nonnegative integers
  we define $M_k = \sum_{i=1}^k m_j$.
\item
$S_0 = \{\nf = \nk: 0<n_1<n_2<n_2< \cdots \}$ denotes the set of all
strictly increasing sequences of natural numbers. $S_0$ is a Baire
space, i.e.\ nonempty open sets are not meager. 
A topological base of open sets is given by all cylinder
sets $[n'_1,n'_2,\ldots,n'_k]$ containing those $\nf \in S_0$
with $n_i=n'_i$ for~$i=1,\ldots,k$.
\item 
$S(\If,\mf)$ defines the closed (and hence Baire) subspace of~$S_0$
containing those $\nf$ having with each $I_j$ exactly $m_j$ members in common.
\item
Admissible $(\If,\mf)$ and~$\qf$:
Given $\If = (I_j)_{j \in \N}$ and $\mf = (m_j)_{j \in \N}$ 
with $m_j \le b_j = |I_j|$, we consider the sequence
$\qf$ of ratios $q_j = \frac{m_j}{b_j} \in [0,1]$ and, for each $N \in \N$,
the measure
\[
\pi_{\If,\mf,N} = M_N^{-1} \sum_{j=1}^N m_j \delta_{q_j}
\]
on~$[0,1]$. $(\If,\mf)$ is called admissible if the further conditions
$\lim_{j \to \infty} b_j = \infty$ and $\limn \frac{m_n}{M_n}=0$ are satisfied.
We only consider admissible~$(\If,\mf)$.
\item
$P(\If,\mf) = A((\pi_{\If,\mf,N})_{N \in \N})$, 
the set of accumulation measures of the~$\pi_{\If,\mf,N}$, $N\in \N$.
\item
For $\pi \in \PP$ let $F_{\pi}: [0,1] \to [0,1]$ be defined by
\[
F_{\pi}(t_0) = \pi([0,t_0]) + t_0 \int_{(t_0,1]} \frac{d\pi(t)}{t}
\]
for $0 \le t_0 \le 1$, in particular $F_{\pi}(0) = \pi(\{0\})$.
\item
For each $\pi \in \PP$ we define the set 
\[
M(\pi,\lambda) = \{\mu \in \MM(X): \ \mu \le 
   F_{\pi} \circ \lambda\}.
\]
\item
$M(\lambda,\If,\mf) = \bigcup_{\pi \in P(\If,\mf)} M(\pi,\lambda)$.
\end{itemize}

Given~$\If$, $\mf$ and~$\xf$,  we are interested in the distribution
behavior of a subsequence $\xf \nf = (x_{n_k})_{k \in \N}$
for typical  $\nf\in S(\If,\mf)$. Theorem \ref{Trw1} shows
that limit measures of such sequences cannot be too far from~$\lambda$,
where the precise statement, of course, depends on the parameters
$\If$ and~$\mf$. Theorem \ref{Trw2} shows that 
everything which might happen, happens typically in the Baire sense, i.e.\
all measures not excluded by Theorem \ref{Trw1} are limit measures of
a generic subsequence $\xf \nf$, $\nf \in S(\If,\mf)$.

Now we are ready to state our results:

\begin{Theo}\label{Trw1}
Let $\xf \approx \lambda$. 
Then $M(\xf \nf) \tm M(\lambda,\If,\mf)$ for all $\nf \in S(\If,\mf)$.
\end{Theo}

\begin{Theo}\label{Trw2}
Suppose
 $\xf \approx \lambda$ with~$\supp(\lambda)=X$.
Then $M(\xf \nf) = M(\lambda,\If,\mf)$ for most $\nf \in S(\If,\mf)$,
i.e.\ the exceptional set of those $\nf$ with
$M(\xf \nf) \neq M(\lambda,\If,\mf)$ is meager in~$S(\If,\mf)$.
\end{Theo}

\subsection{Preliminaries}

Not surprisingly, a rigorous proof of Theorems \ref{Trw1}  and \ref{Trw2}
requires a lot of lemmata and combinatorial technicalities. 
Several of these auxiliary results are collected in:

\begin{Lemma}\label{Lrw.aux}
\begin{enumerate}
\item Given at most countably many $\mu_i \in \MM(X)$, there is an open
  basis of~$X$ contained in $\bigcap_{n \in \N}\CC(\mu_n)$.
\item Let $\mu, \mu_1,\mu_2,\ldots \in \MM(X)$ and $V$ a 
  neighborhood of~$\mu$. Then there is an~$\e>0$ and a finite partition 
  $X = A_1 \cup \cdots \cup A_s$,
  $A_i \in \bigcap_{n \in \N} \CC(\mu_n)$ such that, for all $\nu \in \MM(X)$,
  $|\nu(A_i)-\mu(A_i)|<\e$ for~$i=1,\ldots,s$, implies $\nu \in V$.
\item For any $\pi \in \PP$,
  \[
  F_{\pi}(t_0) = \pi([0,t_0]) + t_0 \int_{(t_0,1]} \frac{d\pi(t)}{t}
  \]
  defines a function $F_{\pi}: [0,1] \to [0,1]$ which is monotonically 
  nondecreasing, continuous and concave.
\item Assume $\lim_n \pi_n = \pi$. Then $\limn F_{\pi_n} = F_{\pi}$ uniformly.
\item If $\mu(A) \le F_{\pi}(\lambda(A))$ for all $A \in \CC(\lambda,\mu)$,
  then the same inequality holds for all Borel sets~$A$, i.e.\
  $\mu \le F_{\pi} \circ \lambda$.
\end{enumerate}
\end{Lemma}

\begin{proof}
\begin{enumerate}
\item Standard.

\item Standard.

\item
Note that, for any $\pi \in \PP$ and fixed $t_0 \in [0,1]$,
$F_{\pi}(t_0) = 
% $ can be written as $
\int_0^1 f_{t_0}(t) d\pi(t)$
with the continuous function $f_{t_0}$ defined by $f_{t_0}(t)=1$
 for $0 \le t \le t_0$
and $f_{t_0}(t) = \frac{t_0}{t}$ for $t_0 < t \le 1$. This shows
$0 \le F_{\pi}(t_0) \le 1$ and that $t_0 \le t_1$ implies $f_{t_0} \le f_{t_1}$
and hence $F_{\pi}(t_0) \le F_{\pi}(t_1)$. Thus
$F_{\pi}: [0,1] \to [0,1]$ is monotonic. Continuity at $0$ follows
from monotonic convergence:
\[
\lim_{t_0 \to 0} F_{\pi}(t_0) = 
   \lim_{t_0 \to 0} \int_{[0,1]} f_{t_0}(t) d\pi(t) =
   \int_{[0,1]} f_0(t) d\pi(t) = F_{\pi}(0).
\]
For other points $0< t_0 < t_1 $ observe that
\[
F_{\pi}(t_1) - F_{\pi}(t_0) = A - B + C
\]
with $A= \pi((t_0,t_1])$, $B= t_0 \int_{(t_0,t_1]} \frac{d\pi(t)}t$
and $C = (t_1-t_0) \int_{(t_1,1]} \frac{d\pi(t)}t$.
For fixed $t_0$ and $t_1 \to t_0$ all three values tend to 0, while
for fixed $t_1$ and $t_0 \to t_1$ both $A$ and $B$ tend to $\pi(\{t_1\})$
and $C$ tends to 0. This shows that $F_{\pi}$ is continuous on the
whole interval~$[0,1]$.

In order to see that $F_{\pi}$ is concave we introduce for
$0<t_0<t_1<t_2 \le 1$ the abbreviations
$A_0=(0,t_0]$, $A_1=(t_0,t_1]$, $A_2=(t_1,t_2]$, $A_3=(t_2,1]$,
$p_i = \pi(A_i)$, and $c_i = \int_{A_i} \frac{d\pi(t)}t$.
It suffices to show that
$2F_{\pi}(t_1) \ge F_{\pi}(t_0) + F_{\pi}(t_2)$ for 
$t_0 = t_1 - \e$ and $t_2 = t_1 + \e$. In this case we have \\
$F_{\pi}(t_0) = p_0 + (t_1-\e)(c_1+c_2+c_3)$,\\
$F_{\pi}(t_1) = p_0 + p_1 + t_1(c_2+c_3)$, and\\
$F_{\pi}(t_2) = p_0 + p_1 + p_2 + (t_1+\e)c_3$.\\
Thus the above inequality reduces to 
$p_1 + t_1c_2 + \e(c_1+c_2) \ge p_2 + t_1c_1$, which follows from
$p_1 \ge (t_1-\e)c_1$ and $p_2 \le (t_1+\e)c_2$.

\item
Using the function $f_{t_0}$ from part (3)
one gets pointwise convergence
\[
\limn F_{\pi_n}(t_0) = \limn \int_{[0,1]} f_{t_0} d\pi_n(t) = F_{\pi}(t_0)
\]
immediately from the definition
of the convergence $\limn \pi_n = \pi$ of measures. 
The uniformity in $t_0 \in [0,1]$ finally follows
from a standard argument on the convergence of monotonic functions.

\item Standard.
%Let $B$ be an arbitrary Borel set and assume, by contradiction,
%$\e = \mu(B) - F_{\pi}\lambda(B) > 0$. By Proposition \ref{Pfpi}
%$F_{\pi}$ is continuous, i.e.\ there is a $\delta>0$ such that
%$\lambda(B) \le t < \lambda(B)+\delta$ implies 
%$F_{\pi}(t) < F_{\pi}\lambda(B) + \frac{\e}2$. By outer regularity
%of $\lambda$ there is an open set $V \tm X$ containing $B$
%with $\lambda(V) < \lambda(B)+\delta$, hence 
%$F_{\pi}\lambda(V) < F_{\pi}\lambda(B) + \frac{\e}2$. By inner
%regularity of~$\mu$ there is a compact $K \tm V$ such that
%$\mu(V) < \mu(K) + \frac{\e}2$. For each $x \in K \tm V$ there is
%an $\e_x>0$ such that the ball $B(x,\e_x)$ with center $x$ and
%radius $\e_x$ is a $\lambda$-$\mu$-continuity set contained in~$V$
%(Proposition \ref{Pstetigkeitsmengen}).
%By compactness of~$K$ there are finitely many $x_1,\ldots,x_n \in K$
%such that $K \tm A = \bigcup_{i=1}^n B(x_i,\e_{x_i})$. Since
%$A \in \CC(\lambda,\mu)$ with $K \tm A \tm V$ we conclude 
%\[
%F_{\pi}\lambda(A) \le F_{\pi}\lambda(V) < F_{\pi}\lambda(B) + 
%  \frac{\e}2 = \mu(B) - \frac{\e}2 \le \mu(V) - \frac{\e}2 < 
%  \mu(K) \le \mu(A),
%\]
%contradicting the assumption of the lemma.
\end{enumerate}
\end{proof}

\subsection{Proof of Theorem \ref{Trw1}} \label{SSrw1.proof}
 % begin{proof}[Proof of Theorem \ref{Trweins}]

\begin{Lemma}\label{Lungleichung}
Suppose that $X$ is finite, $A \tm X$, $\nf \in S(\If,\mf)$,
$(\If,\mf)$ admissible, $\lim_{k \to \infty} \mu_{\xf \nf,M_{N_k}} = \mu$,
$\lim_{k \to \infty} \pi_{\If,\mf,N_k} = \pi$ for some 
$N_1 < N_2 < \cdots \in \N$, and $\lambda(A) \le t_0$. 
Then $\mu(A) \le F_{\pi}(t_0)$.
\end{Lemma}
\begin{proof}
By continuity of~$F_{\pi}$ (Lemma \ref{Lrw.aux}(3)) 
it suffices to prove the statement for $\lambda(A) = t_0>0$.
Fix any~$\e>0$. Since $\xf \approx \lambda$ and 
$\lim_{j \to \infty} b_j = \infty$ there is a $j(\e) \in \N$ such that, 
letting  $c_j := |\{n_k\in I_j: x_{n_k} \in A\}|$, we have 
\[
b_j(t_0-\e) \le c_j \le b_j(t_0+\e)
\]
for all $j > j(\e)$. For fixed $k \in \N$ define
\[
\begin{array}{lll}
J_1 &=& \{1,\ldots,j(\e)\}, \\
J_2 &=& \{j:\ j(\e) < j \le M_{N_k}:\ q_j = \frac{m_j}{b_j}
  \le \lambda(A) = t_0\},\  \mbox{and}\\
J_3 &=& \{j:\ j(\e) < j \le M_{N_k}:\ q_j = \frac{m_j}{b_j} 
  > \lambda(A) =t_0\}.
\end{array}
\]
 % To compute $\mu(A) = \lim_{k \to \infty} \mu_{\xf,M_{N_k}}(A)$ 
We are
going to estimate 
$C=C_1+C_2+C_3$, $C_i = \sum_{j \in J_i} c_j$,
$i=1,2,3$.  % , $c_j = |\{n_k \in I_j:\ x_{n_k} \in A\}|$. 
Abbreviate $\pi_{\If,\mf,N_k}$ by~$\pi'_k$.  Now  
$C_1 \le M_{j(\e)}$ is a constant not depending on~$k$,
\[
C_2 \le \sum_{j \in J_2} m_j \le M_{n_k} \pi'_k([0,t_0]),
\]
and, using $b_j \le \frac{m_j}{t_0}$ for $j \in C_3$,
\[
\begin{array}{lll}
C_3 &\le& \sum_{j \in J_3} b_j(t_0 + \e) \le t_0 
     \sum_{j \in J_3} \frac{m_j}{q_j} + 
     \frac{\e}{t_0} \sum_{j \in J_3}m_j \le\\
  &\le& t_0M_{N_k} \int_{(t_0,1]} \frac{d\pi'_k(t)}{t} + \frac{\e}{t_0}M_{N_k}.
\end{array}
\]
Hence
\[
\lim_{k \to \infty} \mu_{\xf,M_{n_k}}(A) \le \lim_{k \to \infty} 
     \frac{M_{j(\e)}}{M_{n_k}} + \pi'_k([0,t_0]) + t_0 
     \int_{(t_0,t]} \frac{d\pi'_k(t)}t
     + \frac{\e}{t_0} = F_{\pi'_k}(t_0) + \frac{\e}{t_0}.
\]
Since this holds for all~$\e>0$, Lemma \ref{Lrw.aux}(4) proves the assertion.
\end{proof}

\begin{Lemma}\label{Lrw1endlich}
Theorem \ref{Trw1} holds whenever $X$ is finite.
\end{Lemma}
\begin{proof}
Let $\mu \in M(\xf \nf)$, $\nf \in S(\If,\mf)$. This means that some
subsequence of the $\mu_{\xf \nf,N}$, $N \in \N$, converges to~$\mu$.
Since $m_k = o(M_k)$ ($k \to \infty$) and since $\PP$ is compact
we may assume $\lim_{k \to \infty} \mu_{\xf \nf,M_{N_k}} = \mu$
and, if necessary again by taking
an appropriate subsequence, $\lim_{k \to \infty} \pi_{\If,\mf,N_k} = \pi$
for some $\pi \in P(\If,\mf)$. Since $X$ is finite, this implies
$\mu(A) \le F_{\pi}\lambda(A)$
for all $A \tm X$. Thus $\mu \in M(\pi,\lambda) \tm M(\lambda,\If,\mf)$.
\end{proof}

\begin{proof}[Proof of Theorem \ref{Trw1}] Let $\mu \in M(\xf \nf)$,
$\nf \in S(\If,\mf)$ and $A \in \CC(\lambda,\mu)$. Consider the
finite space $X' = \{0,1\}$ with measures $\mu'(\{1\}) = \mu(A)$,
$\lambda'(\{1\}) = \lambda(A)$ and the sequence $\xf' = (x'_n)_{n \in \N}$
with $x'_n=1$ iff $x_n \in A$. Then we are in the situation of
Lemma \ref{Lrw1endlich} to conclude $\mu' \in M(\pi,\lambda')$
for some $\pi \in P(\If,\nf)$, i.e.\ 
\[
\mu(A) = \mu'(\{1\}) \le F_{\pi}\lambda'(\{1\}) = F_{\pi}\lambda(A).
\]
Since this holds for arbitrary $A \in \CC(\lambda,\mu)$,
Lemma \ref{Lrw.aux}(5) yields $\mu \in M(\pi,\lambda) \tm M(\lambda,\If,\mf)$.
\end{proof}

\subsection{Proof of Theorem \ref{Trw2}}\label{SSrw2.proof}

\begin{Lemma}\label{Lfortsetzungsfolge}
Assume that $X = \{1,2,\ldots,s\}$ is finite, 
$\lim_{k \to \infty} \pi_{\If,\mf,N_k} = \pi$,
$\e > 0$, $\mu \le F_{\pi}\lambda$, $\nf' \in S(\If,\mf)$,
$\xf \approx \lambda$, $\mbox{supp}(\lambda)=X$ and $j_0 \in \N$. 
Then there exist $j_1 \ge j_0$, $\nf \in S(\If,\mf)$, with $n'_k=n_k$
for all $n_k \in \bigcup_{j=1}^{j_0}I_j$ and 
$|\mu(i) - \mu_{\xf \nf,M}(i)| < \e$ for $M = \sum_{j=1}^{j_1}m_j$ and
all $i \in X$.
\end{Lemma}
\begin{proof}
Let $n_1 < n_2 < \cdots < n_{M_0}$ be given such that
$M_0 = \sum_{j=1}^{j_0}m_j$ and $|I'_j| = m_j$ for $j=1,\ldots,j_0$
and $I'_j = I_j \cap \{n_k:\ 1 \le k \le M_0\}$. We assume
$\lim_{k \to \infty} \pi_{\If,\mf,N_k} = \pi$ for
$N_1 < N_2 < \cdots$, furthermore $\mu \le F_{\pi}\lambda$ and~$\e>0$.
We have to find a number $j_1=N_k \ge j_0$ and an extension 
$n_1 < n_2 < \cdots < n_{M_0} < \cdots < n_M$, $M = \sum_{j=1}^{j_1}m_j$
such that, putting $I'_j = I_j \cap \{n_k:\ k \in \N\}$, we have 
$|I'_j| = m_j$ for all $j=1,\ldots,j_1$ (we call such an
$\nf = (n_1,\ldots,n_M)$ admissible) and $|d_i| < \e$ for all
$d_i = \mu(i) - \mu_{\xf \nf,M}(i)$, $i \in X$. 

To do this let, w.l.o.g., 
$\lambda(1) = \min_{i \in X} \lambda(i)$ which is positive since $X$ is finite
and $\mbox{supp}(\lambda)=X$. Define $c_j^{(i)} = |\{n \in I_j:\ x_n = i\}|$.
Since $\xf \approx \lambda$ and 
$\lim_{j \to \infty} b_j = \infty$ there is some $j' \ge j_0$ such that
$|\frac{c_j^{(i)}}{b_j} - \lambda(i)| < \frac{\lambda(1) \e}{3s^2}$ for 
all $j \ge j'$. Choose $j_1 =N_k > j'$ such that $M'/M < \e/3s$  
($M = \sum_{j=1}^{j_1}m_j$, 
$M' = \sum_{j=1}^{j'}m_j$), $\frac 1M < \frac{\e}{2s^2}$ and
$F_{\pi_{\If,\mf,N_k}} > F_{\pi} - \frac{\e}3$ (Lemma \ref{Lrw.aux}(4)). 
Let now, for our given admissible~$\nf$, 
$A(\nf,i) =  \mu_{\xf \nf,M}(i) = |\{k \le M:\ x_{n_k}=i\}|$. 
Rearrange the $d_i$ in such a way that 
$d_{i_1} \ge d_{i_2} \ge \cdots \ge d_{i_s}$. Since the set of admissible
$\nf$ is finite, there is a nonempty set of admissible $\nf$
for which  $D(\nf) = \sum_{i=1}^s |d_i|$ takes a minimal value, say~$D_0$. 
Among 
these $\nf$ choose one which leads to the $s$-tuple $(d_{i_1},\ldots,d_{i_s})$ 
which is minimal with respect to the lexicographic ordering.
Everything we have to show is $d_{i_1} < \e/s$, since then 
$\sum_{i \in X}d_i = 0$ implies $\max_{i \in X} |d_i| = 
\max \{d_{i_1},-d_{i_s}\} \le \e$. 

Assume therefore, by contradiction, $|d_{i_1}| \ge \e/s$. We treat only the
case $d_{i_1} \ge \e/s > 0$, since $d_{i_1} \le -\e/s < 0$ is similar.
$\sum_{i=1}^s d_i = 0$ and $d_{i_1} \ge \cdots \ge d_{i_s}$
implies $d_{i_s} \le 0$. It follows that there is some $r$ 
such that $d_r - d_{r+1} > \frac{\e}{s^2}$ and $d_r > \frac{\e}{s^2}$.
Let $Y = \{i_1,\ldots,i_r\}$. For $j = j_0+1,\ldots,j_1$ we claim that \\
(i) $x_n \in Y$ for all $n \in I'_j$ whenever 
$\sum_{i \in Y}c_j^{(i)} \ge m_j$;\\
(ii) $x_n \notin Y$ for all $n \in I_j \setminus I'_j$ whenever
$\sum_{i \in Y}c_j^{(i)} < m_j$.\\
To see this, note that, if (i) failed by some $x_n \notin Y$, $n \in I'_j$,
$\sum_{i \in Y}c_j^{(i)} \ge m_j$, we could replace $n$ in~$I'_j$ by some
$n' \in I_j$ with $x_{n'} \in Y$ to get a contradiction to the
extremal choice of~$\nf$. A similar argument shows (ii).

Let now $A_j = |\{n_k \in I_j:\ x_{n_k} \in Y\}|$. Then,  
for $j' \le j \le j_1$, (i) and (ii) together with the extremal choice 
of $\nf$ guarantee the following two implications:\\
$q_j \le \lambda(Y)$ implies $A_j \ge m_j(1- \frac{\lambda_1 \e}{3s})
   \ge m_j (1 - \frac{\e}{3s})$ and \\
$q_j > \lambda(Y)$ implies 
 $A_j \ge (\lambda(Y) - \frac{\lambda_1 \e}{3s})b_j$.\\
Since $d_i \ge 0$ for all $i \in Y$ and $i_1 \in Y$,
$d_{i_1} \ge \e/s$ implies
\[
\mu(Y) \ge  \mu_{\xf \nf,M}(Y) + \frac{\e}s.
\]
In
\[
\mu_{\xf \nf,M}(Y) = \frac 1M \sum_{j \le j_1}A_j 
\]
we split the sum $\sum_{j \le j_1} A_j$ into three sums~$S_i$,
 $i=0,1,2$, where
the summation runs over all $j \in J_i$. Here $J_0$ contains all
$j < j'$, $J_1$ all $j$ with $j' \le j \le j_1$
and $q_j \le \lambda(Y)$, $J_2$ all $j$ with $j' \le j \le j_1$
and $q_j > \lambda(Y)$. By using our lower bounds for
$A_j$ if $j \in J_1$ resp.\ if $j \in J_2$ we get
\[
\mu(Y) \ge \frac 1M \left(
	   \sum_{j \in J_1} m_j(1-\frac{\e}{3s}) +
	   \sum_{j \in J_2} (\lambda(Y) -
	      \frac{\lambda_1 \e}{3s})b_j \right) + \frac{\e}s.
\]
Write now $J'_i$ for the set of all $j \le j_1$ (including those $< j'$)
with $q_j \le \lambda(Y)$ (if $i=1$)
resp.\ $q_j > \lambda(Y)$ (if $i=2$), $S'_i = \sum_{j \in J_i} m_j$, $i=1,2$.
By separating positive and negative terms in the above inequality
and by using $M'/M < \frac{\e}{3s}$ and, if $q_j \le \lambda(Y)$,
$\lambda(Y) b_j \ge q_jb_j = m_j$ and $b_j < \frac{m_j}{\lambda(Y)} \le
\frac{m_j}{\lambda_1}$, we can continue our estimation with
\[
\mu(Y) \ge \frac 1M(S'_1+S'_2) - \frac{\e}{3sM}(S'_1+S'_2) + 
   \frac{\e}s - \frac{\e}{3s}.
\]
Note that $\frac 1M(S'_1+S'_2) = F_{\pi_{\If,\mf,N_k}}(\lambda(Y))
> F_{\pi}(\lambda(Y)) - \frac{\e}{3s}$ and $S'_1+S'_2 \le \sum_j m_j = M$
to finally obtain
\[
\mu(Y) > F_{\pi}(\lambda(Y)),
\] 
contradicting $\mu \le F_{\pi}\lambda$.
\end{proof}

\begin{Lemma}\label{Lrw2endlich}
%Lemma 4
Theorem \ref{Trw2} holds whenever $X$ is finite.
\end{Lemma}
\begin{proof}
Recall that for $n_1 < n_2 < \cdots < n_l \in \N$ 
the symbol $[n_1,\ldots,n_l]$ denotes
the set of all $\nf$ with this initial part.
The family  of these cylinder sets forms
an open base for the topology on~$S_0$
and their intersections with $S(\If,\mf)$ form an
open base for the subspace~$S(\If,\mf)$. 
Thus the following statement is just a reformulation of
Lemma \ref{Lfortsetzungsfolge}:
\begin{quote}
For all nonempty open sets $O \tm S(\If,\mf)$,
$\mu \le F_{\pi}\lambda$, $\pi \in P(\If,\mf)$, neighborhoods $V$ of~$\mu$
and $l_0 \in \N$ there is a nonempty open set $U \tm O$ and some 
$l \ge l_0$ such that $\nf \in U$ implies $\mu_{\xf \nf,l} \in V$.
\end{quote}
This implies that the set $S(V,l_0,\mu)$ of all $\nf \in S(\If,\mf)$
such that $\mu_{\xf \nf,l} \in V$ for some $l \ge l_0$ is a residual
subset of~$S(\If,\mf)$. Since $\mu$ has a countable
neighborhood base consisting of some~$V_1,V_2,\ldots$, the intersection
$S(\mu)$ of all~$S(V_n,l_0,\mu)$, $n,l_0 \in \N$, is residual in~$S(\If,\mf)$.
Furthermore $\mu \in M(\xf \nf)$ for all $\nf \in S(\mu)$. 

Finally take a sequence $\mu_1,\mu_2,\ldots$ 
which is dense in~$M(\lambda,\If,\mf)$.
Then $S' = \bigcap_{i=1}^{\infty} S(\mu_i)$ is residual and $\nf \in S'$
implies $\mu_i \in M(\xf \nf)$ for all $i$, hence
\[
M(\lambda,\If,\mf) = \overline{\{\mu_i:\ i \in \N\}} \tm M(\xf \nf).
\]
This together with Theorem \ref{Trw1} proves the lemma.
\end{proof}

\begin{proof}[Proof of Theorem \ref{Trw2}]
First fix any
$\mu \in M(\lambda,\If,\mf)$, say $\mu \in M(\pi,\lambda)$ for some 
fixed $\pi \in P(\If,\mf)$. For any given neighborhood $V$ of~$\mu$ 
there is a partition $X=A_1 \cup \cdots \cup A_s$, $A_i \in \CC(\lambda,\mu)$
(cf.\ Lemma \ref{Lrw.aux}, parts(1) and (2)), such that
$\nu(A_i)=\mu(A_i)$ implies $\nu \in V$. Similar to the proof of
Theorem \ref{Trw1} we can here apply the finite case
(Lemma \ref{Lrw2endlich}) to the induced structure 
and see that the set $S(V)$ of
all $\nf \in S(\If,\mf)$ with $M(\xf \nf) \cap V \neq \emptyset$
is residual. By considering the residual set
$S(\mu) = \bigcap_{n=1}^{\infty} S(V_n)$, where $V_1,V_2,\ldots$
form an neighborhood base of~$\mu$, and using that all $M(\xf \nf)$
are closed, one sees that $\mu \in M(\xf \nf)$ for all $\nf \in S(\mu)$.
Now the same argument as in the proof of
Lemma \ref{Lrw2endlich} implies the theorem.
\end{proof}

% \newpage

\end{document}